\documentclass{aptpub}
\usepackage{graphicx}
\usepackage{pstricks,pst-fill,pst-node,pst-text,pst-plot,pst-coil}
\usepackage{wasysym}
\usepackage{mathrsfs}
\usepackage{marvosym}

\authornames{Pierre Calka, Julien Michel, Katy Paroux} 
\shorttitle{Refined convergence for the Boolean model} 

\newcommand{\cqfd}{\par~\hfill$\blacksquare $}
\newcommand{\RR}{{\mathbb R}}
\newcommand{\xxx}{{\mathbf{X}}}

\newcommand{\zzz}{{\mathbf{Z}}}
\newcommand{\converge}[2]{\mathop{\longrightarrow}\limits^{\text{#1}}_{#2}}
\newcommand{\un}[1]{{\text{\bf 1}}_{#1}}
\newcommand{\dist}{\mathop{\rm dist}\nolimits}
\renewcommand{\P}{\mathbf P}
\renewcommand{\E}{\mathbf E}
\newcommand{\acos}{\mathop{\rm acos}\nolimits}

\def\d{\mathrm d}

\newrgbcolor{grey}{0.8 0.8 0.8}

\begin{document}

\title{Refined convergence for the Boolean model} 

\authorone[MAP 5, Universit\'e Paris 5]{Pierre Calka}
\authortwo[UMPA, ENS Lyon]{Julien Michel} 
\authorthree[Laboratoire de Math\'ematiques, Universit\'e de Franche Comt\'e and INRIA 
Rennes - Bretagne Atlantique]{Katy Paroux}

\addressone{MAP5, U.F.R. de Math\'ematiques et Informatique
Universit\'e Paris Descartes
45, rue des Saints-P\`eres
75270  Paris Cedex 06
France, research partially supported by the French ANR project "mipomodim" No. ANR-05-BLAN-0017.
}
\addresstwo{Unit\'e de Math\'ematiques Pures et Appliqu\'ees UMR 5669, ENS Lyon, 
46 all\'ee d'Italie, F-69364 Lyon Cedex 07, France.} 
\addressthree{Laboratoire de Math\'ematiques de Besan\c con, UMR 6623, F-25030 
Besan\c con Cedex, France, and
INRIA,
centre de Rennes - Bretagne Atlantique
Campus de Beaulieu,
F-35 042 Rennes Cedex,
France.}

\begin{abstract}
In a previous work \cite{shell} two of the authors proposed a new proof of a well known convergence result
for the scaled elementary connected vacant component in the high intensity Boolean model towards the
Crofton cell of the Poisson hyperplane process (see e.g.
\cite{hall85}). In this paper, we investigate
the second-order term in this convergence when the two-dimensional Boolean model and the
Poisson line process are coupled on the same probability space. We consider the particular case 
where the grains are discs with random radii.
A precise coupling between the Boolean model and the Poisson line process is first
established. A result of directional convergence in distribution for the
difference of the two sets involved is then derived. Finally we show the convergence of this
directional approximate defect process.
\end{abstract}



\section{Introduction and notations} 

Since the first result of P. Hall \cite{hall85,hall88} and its
generalizations in \cite{molchanov,these}, the scaled
vacancy of the Boolean model is known to converge in some sense to its counterpart in
the Poisson hyperplane process. In a previous paper \cite{shell} two of the authors gave 
another proof of this convergence result for the local occupation laws of a Boolean
shell model in terms of Hausdorff distance. This convergence appears as a first order 
result, expressed in terms of weak convergence. Our aim in this paper is to give two 
generalisations of this result. We extend first the weak convergence to an almost sure
convergence thanks to an adequate coupling between both models, and secondly we
show a second order weak convergence for the difference of both sets, expressed as 
the convergence of a stochastic process in the Skorohod and $L^1$ senses.\\

We shall work in the plane $\RR^{2}$, though some of our results might be stated in higher dimensions: 
let us consider a Boolean model
based on a Poisson point process $\xxx_{\lambda}$ 
with intensity measure $\lambda^2\,\d x$ and generic shape
an open disc centred at $0$ of random radius $R$ such that $\E[R]=1$ and
$\E[R^{2}]<+\infty$.
The law of $R$ will be denoted by $\mu$, and we will assume that there exists $R_{\star}>0$ such
that $\mu(R_{\star},+\infty)=1$.\\
The choice of a random disc enables us to write simply the different couplings and computations presented
below, generic convex smooth shapes could probably be treated in the same way, up to technical
details.\\

The occupied phase of the Boolean
model is denoted by
$${\mathscr O}_{\lambda}=\bigcup_{x\in\xxx_{\lambda}}B(x,R_{x}),$$
where $B(x,r)$ denotes the disc centred at $x$ and of radius $r$ and where
the radii $R_{x}$ for each $x\in\xxx_{\lambda}$ are independent and identically distributed, independent of
$\xxx_{\lambda}$ \cite{molchanov-livre}. 
This process is supposed to leave the point $0$ uncovered, which occurs 
with positive probability
$$\P(0\notin{\mathscr O}_{\lambda})=\exp(-\pi\lambda^2\E[R^{2}]).$$
From now on the Boolean model shall be conditioned by this event.\\
Let $D_{0}^{\lambda}$ denote the (closed) connected component of 
$\RR^{2}\setminus {\mathscr O}_{\lambda}$ containing $0$. The following
asymptotic result for this process (see \cite{hall85,molchanov,these,shell}) may be seen as a consequence
of Steiner's formula \cite{schneider}:
\begin{theorem}\label{cvorigine}
Let $D^\lambda$ be the following compact set:
\begin{itemize}
\item $D^\lambda=\lambda^2D_0^\lambda$ whenever this set is bounded,
\item $D^\lambda=K_0$ a given fixed compact set otherwise.
\end{itemize}
When $\lambda$ tends to infinity, $D^\lambda$ converges in law towards the
Crofton cell ${\mathscr C}$ of
a Poisson line process with intensity measure $\d\rho\,\d\theta$.
\end{theorem}
In \cite{hall85} the convergence was stated for random discs and Hausdorff distance, whereas in \cite{molchanov} it was proved
for generic shapes, using the hit or miss topology for random closed sets. The criterion developed
in \cite{these} gives the convergence for another general class of shapes, whereas in \cite{shell} the convergence is proved
using a convergence result for random shells.\\
The Poisson line process with intensity measure $\d\rho\,\d\theta$ in $\RR^{2}$ is defined as the 
set of (random) lines $D_{\rho,\theta}=\{(r,t)\in{\mathbb R}_{+}\times[0,2\pi)\ :\ r\cos(t-\theta)=\rho\}$, 
where $(\rho,\theta)$ are the points of a Poisson
point process $\Phi$ in ${\mathbb R}_+\times [0,2\pi)$ with intensity measure $\d\rho\,\d\theta$. The Crofton cell ${\mathscr C}$
is defined as the polygon formed
by those lines containing the point $0$ (see \cite{stoyan} or \cite{these} for a survey on
Poisson line tessellations). Numerous distributional results on this model have been
obtained notably by R. E. Miles \cite{miles1,miles2} and G. Matheron
\cite{matheron}.  More recently, central limit theorems have been derived in
\cite{these,katy-tcl} for the two-dimensional case and in \cite{schmidt-tcl} for the
general case. Besides, D. G. Kendall's conjecture on the shape of the Crofton cell when it
is large has been proved in \cite{hug}. Additional distributional and
asymptotic results at large inner radius have also been obtained in \cite{pierre-loi-cotes} and
\cite{pierre-loi-rayon}.\\

We shall first recall in section 2 the asymptotic properties of the outer radius of the
Crofton cell and give some counterpart of those asymptotics for the rescaled outer radius of the
Boolean empty connected component.
Those results are useful for the 
next sections 3, 4 and 5, they also give some
insight on the behaviour of the high intensity Boolean model with respect to the continuous percolation 
problem. It is a natural problem to try to estimate the error 
in theorem \ref{cvorigine}: one possible answer is to give a geometric description of the difference of
those two sets. This description requires to couple the Boolean model with the Poisson line
process. This coupling, which asserts as a consequence the almost sure convergence in Theorem
\ref{cvorigine}, will be described in section 3, and its application to the second order convergence
will be treated in section 4 for directional convergence, and in section 5 for the convergence of the 
{\em rescaled defect process} in the Skorohod and $L^1$ settings.\\

\section{Estimates on the tail probability of the inner and outer radius}

This section contains autonomous results about the inner and outer radii of both the Crofton cell and the 
the empty connected component $\lambda^2 D_{0}^\lambda$. Let us introduce some notations:
\begin{itemize}
\item {\it Crofton cell:\/} the inner radius is denoted by $R_{m}$, the outer radius by $R_{M}$;
\item {\it rescaled Boolean model:\/} the inner radius is denoted by $R_{m}(\lambda)$, the outer radius by $R_{M}(\lambda)$;
\end{itemize}
and they are defined by
\begin{eqnarray*}
R_{m}=\sup\{r>0\ :\ B_{2}(0,r)\subset {\mathscr C}\},&&R_{M}=\inf\{r>0\ :\ {\mathscr C}\subset B_{2}(0,r)\},\\
R_{m}(\lambda)=\sup\{r>0\ :\ B_{2}(0,r)\subset \lambda^2 D_{0}^\lambda\},&&R_{M}(\lambda)=
\inf\{r>0\ :\ {\lambda^2 D_{0}^\lambda}\subset B_{2}(0,r)\},\\
\end{eqnarray*}
The laws of some of those quantities are well known and straightforward to obtain:
\begin{eqnarray*}
\forall r>0,&&\P(R_{m}>r)=\exp(-2\pi r),\\
&&\P(R_{m}(\lambda)>r)=\exp(-(2\pi r+\pi r^2/\lambda^2)),
\end{eqnarray*}
however for the outer radii we only have the following asymptotic result, proved in \cite{pierre-loi-rayon} in the context 
of a study of Kendall's conjecture on the shape of large Poisson polygons:
\begin{theorem}[Theorem 8 in \cite{pierre-loi-rayon}]\label{loi-pierre-rayon}
We have
for all $r>0$
\begin{eqnarray*}
&&2\pi re^{-2r}\left(\cos 1+\frac{e^{-(2\pi\cos 1-1)r}}{2\pi r}\right)\leq
 \P(R_{M}\geq r)\\
 &\leq &2\pi re^{-2r}\left(1-(\pi-2)re^{-2r}+\frac{2}{3}(\pi-3)^{2}r^{2}e^{-4r}
+\frac{e^{-2(\pi-1)r}}{2\pi r}\right).
\end{eqnarray*}
\end{theorem}
Concerning the Boolean model, we prove below the following counterpart:
\begin{theorem}\label{theo-bool-RM}
There exists a constant $S>0$ and constants $K$ and $C$ depending only on $R_{\star}$
such that for all $r\geq 1$ and $\lambda^{2}\geq Sr$ we have
$$\P(R_{M}(\lambda)>r)\leq K
\exp\left(-Cr\right).$$
\end{theorem}
\begin{proof}
The proof of this theorem relies on the following (non-optimal)
reasoning~:
let $N$ be a positive integer greater than $12$, and define the angular sectors $S_{i,N}$ for $i\in\{1,\ldots,N\}$ as
$$S_{i,N}=\left\{(\rho,\theta)\,:\,\rho>0,
\theta+\frac{\pi}{N}\in \left[\frac{(2i-1)\pi}{N},\frac{(2i+1)\pi}{N}\right)\right\},$$
for $i\in\{0,\ldots,N-1\}$. If $\lambda^2 D_{0}^\lambda$ is not included
in $B(0,r)$, then there exists at least one of those sectors such that no disc of the rescaled Boolean model 
contains both points with polar coordinates $(r,(2i-1)\pi/N)$ and 
$(r,(2i+1)\pi/N)$ (this implies that $\lambda$ must satisfy the condition
$\lambda^{2}R\geq r\sin\pi/N$).\\

Let us denote
by $A_{\lambda,r,R}$ the set of the centres $(\rho,\theta)$ of discs of radius $\lambda^2 R$ 
such that this occurs for the sector $S_{0,N}$, we have by invariance under rotations
\begin{eqnarray}
\P(R_{M}(\lambda)\geq r)&\leq& N \exp\left(-
\lambda^{-2}\int_{\RR_{+}\times[0,2\pi)\times\RR_{+}}
\un{A_{\lambda,r,R}}(\rho,\theta)\,\rho\,\d\rho\,\d\theta\,\d\mu(R)
\right),\label{majo-theo-bool-RM}
\end{eqnarray}
and the aim of the computations is to bound the Lebesgue measure $a_{\lambda,r,R}$
of the set $A_{\lambda,r,R}$ from
below. This set is
$$A_{\lambda,r,R}=\left[B\left((r,\pi/N),\lambda^{2}R\right)
\cap B\left((r,-\pi/N),\lambda^{2}R\right)\right]
\setminus B\left(0,\lambda^{2}R\right).$$

The geometry of this set is quite easily described, let us indeed introduce the angle
$\theta_{0}=\acos (r/2\lambda^2R)-\pi/N$, then for $\lambda$ such that
\begin{equation}\label{est-lambda-RM}
\lambda^2R\geq \max\left(\frac{r}{2\cos(\theta_{0}-\pi/N)},\frac{r}{2\cos(\pi/N)}\right),
\end{equation}
one has
$$(\rho,\theta)\in A_{\lambda,r,R}\text{ if and only if } \left\{\begin{array}{l}
-\theta_{0}\leq \theta\leq \theta_{0},\\
\lambda^2R\leq\rho\leq \rho_{e}(\theta),
\end{array}\right.$$
where $\rho_{e}(\theta)$ is given by
$$\rho_{e}(\theta)=
r\cos(|\theta|+\pi/N)+\lambda^2R\sqrt{1-\left(\frac{r\sin(|\theta|+\pi/N)}{\lambda^2R}\right)^2}.$$
Let us now introduce $\theta_{1}=\pi/12$, this angle satisfies $\theta_{1}<\theta_{0}$
for $\lambda$ large enough ($\lambda^2\geq Sr$ where the constant
$S$ is chosen greater than $2$, and depends on $N$ and $R_{\star}$). The computation of
$a_{\lambda,r,R}$ becomes 
\begin{eqnarray*}
a_{\lambda,r,R}&=&\int_{0}^{\theta_0}\left[\rho_{e}(\theta)^2-(\lambda^2R)^2\right]\,d\theta,\\
&\geq&\int_{0}^{\theta_1}\left[\rho_{e}(\theta)^2-(\lambda^2R)^2\right]\,d\theta,\\
&\geq&\int_{0}^{\theta_1}\Bigl[\underbrace{r^2\cos(2\theta+{2\pi}/{N})}_{\geq 0}+\\
&&2r\cos(\theta+\frac{\pi}{N})\underbrace{\sqrt{\lambda^4R^2-
r^2\sin^2(\theta+{\pi}/{N})}}_{\geq (\sqrt{3}/2)\lambda^{2}R\text{ from (\ref{est-lambda-RM})}}
\Bigr]\,d\theta,\\
&\geq&{\frac{3\pi}{24}}\,\lambda^2Rr.
\end{eqnarray*}
Consequently we obtain that there exists a constant $C'>0$ such that
$${a}_{\lambda,r,R}\geq C'\lambda^2R r.$$
Inserting this estimate in inequality (\ref{majo-theo-bool-RM})
completes the proof of theorem \ref{theo-bool-RM}.

\cqfd\end{proof}

\section{Coupling and almost sure convergence}

\begin{figure}[h]
\begin{center}
\begin{picture}(300.00,132.00)
\includegraphics{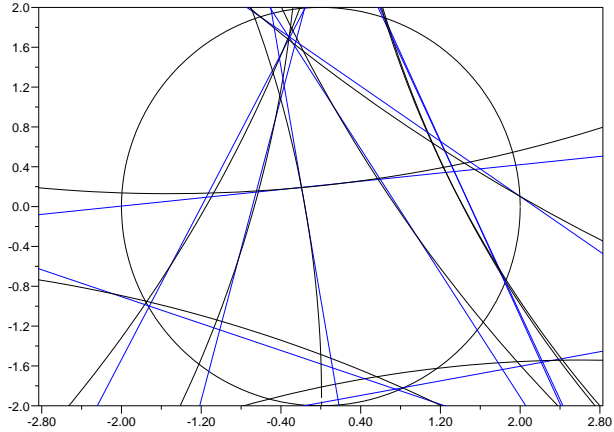}
\end{picture}
\end{center}
\vspace*{-5cm}
\caption{Coupling the Boolean model with the line process (the simulation is exact within the circle).}\label{couplage-dessin}
\end{figure}

Coupling the Boolean model with the Poisson line process is an easy task: indeed as $\lambda$
tends to infinity the rescaled Boolean model {\em looks like} the Poisson line process as one can
see from theorem 1 in \cite{shell}. The formal way to state this as a coupling result
is to introduce a marked Poisson line process which couples both processes:\\
let $\zzz$ be a Poisson point process with intensity measure $\d\rho\,\d\theta\,\d\mu$
on $\RR_{+}\times[0,2\pi)\times\RR_{+}$, and define the function 
$$\psi_{\lambda}(\rho,R)=R\sqrt{1+\frac{2\rho}{\lambda^2R}}.$$
Define the following processes:
\begin{itemize}
\item $\xxx=\bigcup_{(\rho,\theta,R)\in\zzz}\{(\rho,\theta)\}$,
\item $\xxx_{\lambda}^M=\bigcup_{(\rho,\theta,R)\in\zzz}\{(\psi_{\lambda}(\rho,R),\theta,R)\}$,
\end{itemize}
then one has:
\begin{proposition}\label{couplage-LB}
$\xxx$ and $\xxx_{\lambda}^M$ are Poisson point processes with respective intensities 
$\d \rho\,\d\theta$, and $\lambda^2\,\rho\,\d\rho\,\d\theta\,\d\mu$
on respectively ${\mathbb R}_{+}\times[0,2\pi)$ and ${\mathscr B}:=\{(r,t,R)\in{\mathbb R}_{+}\times[0,2\pi)\times
{\mathbb R}_{+}\ :\ r>R\}$. 
\end{proposition}
We may then construct the polar lines at the points of $\xxx$, the Boolean model of discs associated to $\xxx_{\lambda}^M$:
this Boolean model does not cover the origin, rescale this last
Boolean model by an homothetic factor $\lambda^2$, and compare them, this procedure is illustrated in figure \ref{couplage-dessin}.
\vspace*{6cm}

\begin{remark}
Conversely we could have introduced the coupling starting from the points $(\rho,\theta,R)$ 
of a marked Poisson point process $\xxx_{\lambda}^M$
with intensity $\lambda^2\,\rho\,\d\rho\,\d\theta\,\d\mu$ on ${\mathscr B}$, yielding directly the 
Boolean model. In this setting the application 
\begin{equation}\label{formule-couplage-BL}
(\rho,\theta,R)\mapsto \left(\lambda^{2}(\rho-R)+\frac{1}{2R}\lambda^{2}(\rho-R)^{2},\theta\right)
\end{equation}
maps $\xxx_{\lambda}^M$ onto a Poisson point process $\xxx$ with intensity measure $\d\rho\,\d\theta$ on which we may construct the
Poisson line process.\\
\end{remark}

\begin{proof}
Proposition \ref{couplage-LB} is easily proved by the following:
\begin{itemize}
\item $\xxx$ is clearly a Poisson point process with the right intensity measure;
\item and $\xxx_{\lambda}$ is also a Poisson point process, whose intensity measure is the
image of the intensity measure of $\zzz$ by the map
$$(\rho,\theta,R)\mapsto (\psi_{\lambda}(\rho,R),\theta,R),$$
a straightforward computation shows the result.
\end{itemize}
\cqfd\end{proof}

From now on we shall use the coupling induced by $\zzz$, thus the set $D^{\lambda}$ will refer to the rescaled
connected component in this Boolean model (when bounded, $K_{0}$ otherwise), and ${\mathscr C}$
will be the Crofton cell in this line process.\\

This coupling yields the following result on the local accuracy on the approximation of the rescaled Boolean model
by the line process:
\begin{proposition}\label{cv-hausdorff}
For all $M>0$ and $\lambda\geq\lambda_{0}(M)=\sqrt{\frac{M}{R_{\star}}}>0$ one has
$$d_{H}^{M}(\lambda^2 D_{0}^{\lambda},{\mathscr C})
\leq \frac{M'^{2}}{R_{\star}\lambda^2}\text{ a.s.},$$
where
$M'=M+{M^2}/({\lambda^2R_{\star}})$, and $d_{H}^{M}$ denotes  the $M$-Hausdorff distance in $B(0,M)$ defined by
\begin{eqnarray*}
d_{H}^{M}(F,G)=\inf\bigl\{\alpha>0&:&(F\oplus B(0,\alpha))
\cap B(0,M)\supset G\cap B(0,M),\\
&&(G\oplus B(0,\alpha))
\cap B(0,M)\supset F\cap B(0,M)\bigr\},
\end{eqnarray*}
and for any two subsets $A$ and $B$, the set $A\oplus B$ denotes their Minkowski sum: 
$A\oplus B=\{x+y,\,(x,y)\in A\times B\}$.
\end{proposition}
\begin{proof}
Let us denote by $(\rho_{1},\alpha_{1},R_{1}),\ldots,(\rho_{K},\alpha_{K},R_{K})$ the points of
$\zzz$ such that the rescaled discs $B_{i}$ associated to those points intersect $B(0,M)$. Since $\lambda\geq
\lambda_{0}(M)$ no such disc can be included in $B(0,M)$. A straightforward computation
with the help of formula (\ref{formule-couplage-BL})
shows that the associated lines intersect the disc $B(0,M')$, where 
$$M'=M+\frac{M^2}{\lambda^2R_{\star}}.$$
More precisely, as is shown in figure 
\ref{approxdessin}, the $M$-Hausdorff distance between the intersection of the circle $\partial B_{i}$ with 
$B(0,M)$ and the intersection of
tangent line $T_{i}$ with $B(0,M)$ is bounded by the distance between the points $A$ and $C$,
defined as respectively the intersection of $T_{i}$ and $\partial B(0,M)$, and
the point on $\partial B_{i}$ aligned with $A$ and the center of $B_{i}$. We have, if we define
the distance from the origin to $B_{i}$ as 
$u_{i}=\lambda^2R_{i}(\sqrt{1+2\rho_{i}/(\lambda^2R_{i})}-1)\leq M$:
\begin{eqnarray*}
AC&=&\sqrt{\lambda^{4}R_{i}^{2}+M^{2}-u_{i}^{2}}-\lambda^{2}R_{i},\\
&=&\frac{M^{2}-u_{i}^{2}}{\sqrt{\lambda^{4}R_{i}^{2}+M^{2}
-u_{i}^{2}}+\lambda^{2}R_{i}},\\
&\leq&\frac{M^{2}}{2\lambda^{2}R_{i}}.
\end{eqnarray*}

\vspace*{3cm}
\psline[linewidth=1pt,linestyle=dashed,linecolor=black]{->}(3,0)(10,0)
\psarc[linecolor=red,linewidth=1.5pt](20,0){13}{165}{195}
\psline[linecolor=red,linewidth=1.5pt](7,3)(7,-3)
\psline[linecolor=blue,linewidth=1.5pt](7.7,3)(7.7,-3)
\psarc[linecolor=black,linewidth=1pt,showpoints=true](4,0){4}{310}{50}
\pscircle[linecolor=black,fillcolor=black,fillstyle=solid](4,0){0.1}
\rput[b](4,-0.6){$O$}
\rput[b](8,2.0){$D_{i}$}
\rput[b](6.6,1.3){$T_{i}$}
\rput[b](5.7,-3.2){$M$}
\rput[b](7.53,3.5){$\partial B_{i}$}
\rput[b](5.6,-0.7){$u_{i}$}
\psline[linecolor=black,linewidth=0.5pt]{->}(6,-0.3)(7,0)
\pscircle[linecolor=red,fillcolor=red,fillstyle=solid](7,0){0.07}
\rput[b](9.3,-0.7){$\rho_{i}$}
\psline[linecolor=black,linewidth=0.5pt]{->}(8.7,-0.3)(7.7,0)
\pscircle[linecolor=blue,fillcolor=blue,fillstyle=solid](7.7,0){0.07}
\pscircle[linecolor=red,fillcolor=red,fillstyle=solid](7,2.646){0.07}
\rput[b](6.7,2.3){$A$}
\psline[linecolor=red,linestyle=dashed](7,2.646)(20,0)
\psline[linecolor=white,fillcolor=white,fillstyle=solid](20,0)(20,2)(13,2)(13,0)(20,0)
\pscircle[linecolor=red,fillcolor=red,fillstyle=solid](7.28,2.59){0.07}
\rput[b](7.5,2.7){$C$}
\psarc[linecolor=red,linewidth=0.5pt,linestyle=dashed](20,0){13.26}{165}{195}
\rput[b](11.2,1.2){$\to$ to the centre}
\rput[b](11.2,0.8){of the disc}

\vspace*{3cm}
\begin{figure}[h]
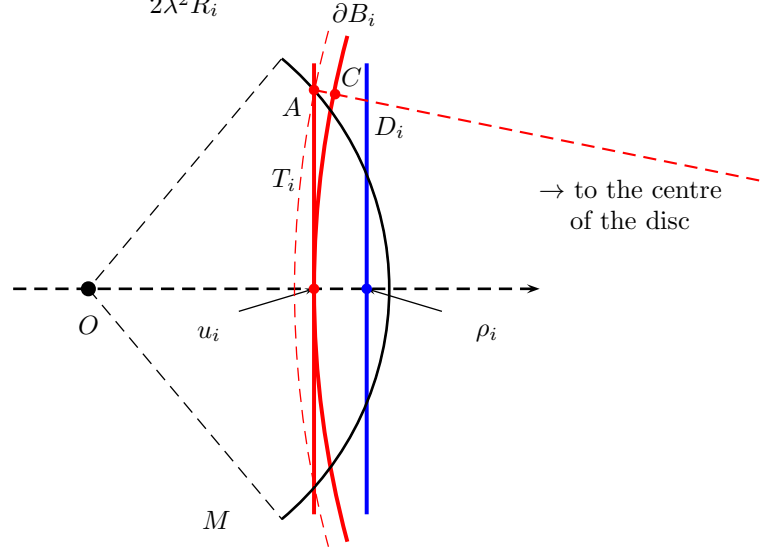

\caption{The circle $\partial B_{i}$ and the tangent line $T_{i}$.}\label{approxdessin}
\end{figure}

Consequently, for any disc $B$ of radius 
$\lambda^2R>\lambda^2R_{\star}$ such that $0\notin B$, and $T$ the tangent to this
disc at its closest point to $0$ one has
\begin{equation}\label{ineg1}
d_{H}^{M}(\partial B,T)\leq \frac{M^{2}}{2\lambda^2R_{\star}}.
\end{equation}
On the other hand, as $u_{i}\leq M$, one has
\begin{eqnarray}
\left|u_{i}-\rho_{i}\right|&=&\rho_{i}\,\left|
\frac{1-\sqrt{1+\frac{2\rho_{i}}{\lambda^2R_{i}}}}{1+\sqrt{1+\frac{2\rho_{i}}{\lambda^2R_{i}}}}\right|\nonumber\\
&\leq&
\frac{\rho_{i}^{2}}{2\lambda^2 R_{i}},\nonumber\\
&\leq&
\frac{M'^{2}}{2\lambda^2 R_{\star}},\label{ineg2}
\end{eqnarray}
so that one obtains the result of proposition \ref{cv-hausdorff} 
by combining the two inequalities (\ref{ineg1},\ref{ineg2}): if we denote 
by $H_{i}$ (resp. $H'_{i}$) the half plane with boundary $D_{i}$ (resp. $T_{i}$)
not containing the origin:
\begin{eqnarray*}
d_{H}^{M}(\lambda^2 D_{0}^{\lambda},C_{0}^{\lambda})
&\leq&\max_{i\in\{1,\ldots,K\}}d_{H}^{M}({}^{\complement}B_{i},{}^{\complement}H_{i}),\\
&\leq&\max_{i\in\{1,\ldots,K\}}d_{H}^{M}({}^{\complement}B_{i},{}^{\complement}H'_{i})\\
&&+\max_{i\in\{1,\ldots,K\}}d_{H}^{M}({}^{\complement}H_{i},{}^{\complement}H'_{i}),\\
&\leq&\frac{M'^{2}}{\lambda^2R_{\star}},
\end{eqnarray*}
where ${}^\complement G$ denotes the complementary set of $G$, and this concludes the proof of proposition \ref{cv-hausdorff}.
\cqfd\end{proof}

From proposition \ref{cv-hausdorff} we may deduce the almost sure convergence in our coupled setting:
\begin{theorem}
Almost surely $D^\lambda$ converges in Hausdorff distance towards ${\mathscr C}$.
\end{theorem}
\begin{proof}
Let us consider the subset $\Omega_{r}$ of those $\omega's$ such that both $R_{M}$ and $R_{M}(\lambda)$ are lesser than $r$,
then the Hausdorff distance between ${\mathscr C}$ and $D^\lambda=\lambda^2D_{0}^\lambda$ 
is lesser than $r'^2/(R_{\star}\lambda^2)$ where 
$r'=r+r^2/(\lambda^2R_{\star})$ for
$\lambda$ large enough thanks to proposition \ref{cv-hausdorff}, and $\P(\Omega_{r})\to 1$ as $r\to+\infty$
thanks to theorems \ref{loi-pierre-rayon} and \ref{theo-bool-RM}, thus the conclusion.
\cqfd\end{proof}
\begin{remark}
The actual speed of convergence shown above could be stated in the following way: let $S_{M}(\lambda)=\max(R_{M},R_{M}(\lambda))$,
then
$$d_{H}({\mathscr C},D^\lambda)\leq \frac{S_{M}'(\lambda)^2}{\lambda^2R_{\star}},$$
where
$S_{M}'(\lambda)=S_{M}(\lambda)+S_{M}(\lambda)^2/(\lambda^2R_{\star})$.
\end{remark}

\section{Convergence of the second order-directional results}

In order to prove a second order convergence result for the empty connected component towards the
Crofton cell, we shall first give some notations and definitions, then we shall state the convergence
results for one, then many directions in a second and third subsections. 

\subsection{Notations}

Recall that for each $\omega\in\Omega$, we denote by ${\mathscr C}(\omega)$ the Crofton cell of the Poisson line process
$\xxx$ induced by $\zzz$, and by
$D^\lambda(\omega)$ the empty connected component of the rescaled coupled Boolean model. We define the following quantities
(almost surely they are all finite random variables):
\begin{itemize}
\item $N_{e}(\omega)$ the number of vertices of ${\mathscr C}(\omega)$, those points are denoted 
anti-clockwise by $V_{1}(\omega),\ldots,V_{N_{e}(\omega)}(\omega)$;
\item  $0\leq \theta_{1}(\omega)<\ldots<\theta_{N_{e}(\omega)}(\omega)<2\pi$ the polar angles of 
those vertices;
\item we take the convention for the edge numbered $i$ of ${\mathscr C}(\omega)$ to 
join vertices $V_{i}$ (included) and $V_{(i\mod N_{e}(\omega))+1}$ (excluded);
\item for each $i\in\{1,\ldots,N_{e}(\omega)\}$, set $(\Upsilon_{i}(\omega),\Theta_{i}(\omega),R_{i}(\omega))$ 
the polar coordinates ({\em angle} and {\em distance}) of the edges of ${\mathscr C}(\omega)$ marked with the 
associated radius of the disc in the coupled Boolean model (from now on we will write $i+1$ for 
$(i\mod N_{e}(\omega))+1$ and $N_{e}$ for $N_{e}(\omega)$, for sake of simplicity).
\end{itemize}
For each $\omega\in\Omega$ and $t\in[0,2\pi)$ we define $\Delta_{t}$ the half-line $\{(r,t)\ :\ r>0\}$ and
\begin{itemize}
\item $\Theta(t,\omega)$ the polar angle of the edge intersecting the half-line $\Delta_{t}$;
\item $\Upsilon(t,\omega)$ the distance from the origin to this edge;
\item $L(t,\omega)=\Upsilon(t,\omega)/\cos(\Theta(t,\omega)-t)$ the distance from the origin to the 
intersection of $\Delta_{t}$ with this edge;
\item $R(t,\omega)$ the radius of the associated disc.
\end{itemize}
All quantities above are well defined on the same set of full probability for each $t$.

\begin{figure}[h]
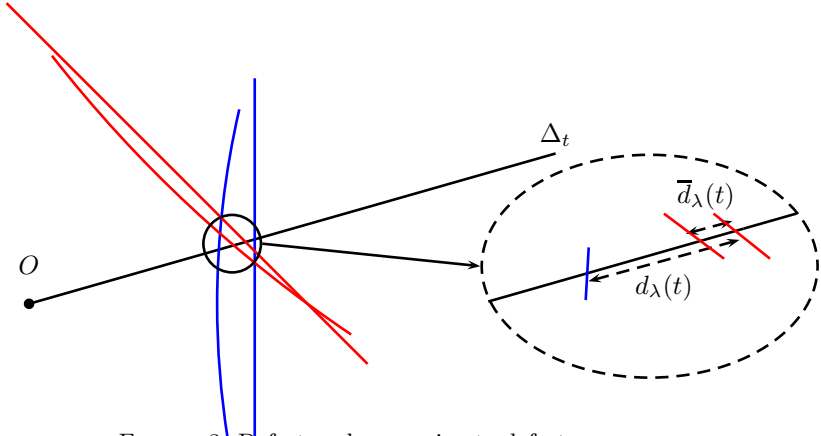

\caption{Defect and approximate defect}\label{a-faire-dessin-1}
\end{figure}

\begin{definition}\label{def-def-app}
For each $\omega\in\Omega$ and $t\in[0,2\pi)$ we define
the {\em defect} at angle $t$ by
$$d_{\lambda}(t,\omega)=\dist(0,\Delta_{t}\cap \partial D^\lambda(\omega))
-\dist(0,\Delta_{t}\cap \partial {\mathscr C}(\omega)),$$
where $\dist$ denotes the Euclidean distance,
and the {\em approximate defect} at angle $t$ by 
\begin{eqnarray*}
\overline{d}_{\lambda}(t,\omega)&=&\dist\left(0,\Delta_{t}\cap 
B\left(\lambda^2\psi_{\lambda}(\Upsilon(t,\omega),\Theta(t,\omega)),\lambda^2R(t,\omega)\right)\right)-\\
&&\dist(0,\Delta_{t}\cap {\mathscr C}(\omega)),
\end{eqnarray*}
when this quantity is well-defined ($\lambda$ large enough).
\end{definition}
We check easily that 
\begin{eqnarray}
\overline{d}_{\lambda}(t,\omega)&=&\lambda^{2}R(t,\omega)\cos(\Theta(t,\omega)-t)
\sqrt{1+\frac{2\Upsilon(t,\omega)}{\lambda^{2}R(t,\omega)}}\nonumber\\
&&-\lambda^{2}R(t,\omega)\sqrt{1-\sin^2{(\Theta(t,\omega)-t)}
\left(1+\frac{2\Upsilon(t,\omega)}{\lambda^{2}R(t,\omega)}\right)}\nonumber\\
&&-\frac{\Upsilon(t,\omega)}{\cos(\Theta(t,\omega)-t)},\label{formule-defaut-app}
\end{eqnarray}
for $\lambda\geq \sqrt{2\Upsilon(t,\omega)\tan^2(\Theta(t,\omega)-t)/R(t,\omega)}$, see figure \ref{a-faire-dessin-1}.

\vspace*{5cm}

\hspace*{-2cm}\psline[linecolor=black,linewidth=1pt]{*-}(4,1)(11,3)
\rput[b](4,1.4){$O$}
\rput[b](11,3.1){$\Delta_{t}$}
\psline[linecolor=blue,linewidth=1pt](7,-0.8)(7,4)
\psarc[linecolor=blue,linewidth=1pt](18,1){11.5}{167}{189}
\psline[linecolor=red,linewidth=1pt](8.5,0.2)(3.7,5)
\psarc[linecolor=red,linewidth=1pt](18,15){17.38}{218}{236}
\pscircle[linecolor=black,linewidth=1pt](6.7,1.8){0.4}
\psellipse[linecolor=black,fillstyle=none,linewidth=1pt,linestyle=dashed](12.25,1.5)(2.25,1.5)
\psline[linecolor=black,linewidth=1pt]{->}(7.1,1.8)(10,1.5)
\psline[linecolor=black,linewidth=1pt](10.105,1.03)(13.5,2)
\psline[linecolor=black,linewidth=1pt](14.2,2.2)(13.5,2)
\psline[linecolor=black,linewidth=1pt,linestyle=dashed]{<->}(11.435,1.30)(13.445,1.85)
\rput[b](12.44,1.05){$d_{\lambda}(t)$}
\psline[linecolor=black,linewidth=1pt,linestyle=dashed]{<->}(12.73,1.932)(13.35,2.1)
\rput[b](13,2.25){$\overline{d}_{\lambda}(t)$}
\psline[linecolor=blue,linewidth=1pt](11.4,1.05)(11.44,1.75)
\psline[linecolor=red,linewidth=1pt](13.85,1.6)(13.1,2.2)
\psline[linecolor=red,linewidth=1pt](13.24,1.6)(12.44,2.2)

\subsection{One directional convergence}

The first result is an almost sure convergence of the defect function in one fixed direction:
\begin{theorem}\label{une-direction}
For all $t\in[0,2\pi)$ one has
$$\lambda^2d_{\lambda}(t,\cdot)\converge{a.s.}{\lambda\to+\infty}Z(t,\cdot),$$
where $Z(t,\cdot)$ is the random variable defined by
$$\forall \omega\in\Omega,\ Z(t,\omega)
=-\frac{L(t,\omega)^{2}}{2}\,\frac{\cos{2}(\Theta(t,\omega)-t)}{R(t,\omega)\cos(\Theta(t,\omega)-t)},$$
and the common law of $(L(t,\cdot),\Theta(t,\cdot),R(t,\cdot))$ is given by
$$\d(L(t,\cdot),\Theta(t,\cdot),R(t,\cdot))(P)
(\ell,\alpha,r)=\pi\exp(-2\pi\ell)\cos(\alpha-t)
\un{\alpha\in(t-\pi/2,t+\pi/2)}\,\d\ell\,\d\alpha\,\d\mu(r).$$
\end{theorem}
\begin{proof}
The proof of this result proceeds in two steps:
\begin{itemize}
\item restrict the probability space to those events such that for $\lambda$ large enough the defect is equal to the approximate
defect;
\item show that those events cover almost surely $\Omega$.
\end{itemize}
\subsubsection*{Step 1: restricted events}
Let us consider $\delta>0$, $\epsilon>0$, $r>0$ and $s>0$, and consider the subset 
$\Omega_{\delta,\epsilon,r,s}$ of all $\omega\in\Omega$ such that
\begin{itemize}
\item $B(0,r)\subset{\mathscr C}(\omega)\subset B(0,s)$;
\item for each $u\in(t-\delta,t+\delta)$, the intersection of ${\mathscr C}(\omega)$ with $\Delta_{u}$ is on the same edge of 
${\mathscr C}(\omega)$;
\item ${\mathscr C}(\omega)\oplus B(0,\epsilon)$ is not intersected by other lines of the Poisson line process than those
on the boundary of ${\mathscr C}(\omega)$.
\end{itemize}
It is quite obvious for geometrical arguments that if $\lambda$ is large enough, in direction $t$ the defect will be exactly 
equal to the approximate defect, as in the disc of radius $s+\epsilon$ the Hausdorff distance between circles and lines
gets smaller as $\lambda$ increases, and thus in direction $t$ the first intersecting line corresponds to the first intersecting disc.
There remains to compute the exact asymptotics of the approximate defect, this is done in the following way, where we restrict ourselves
thanks to invariance under rotations, to the angle $t=0$.\\
One has on the one hand the following well-known classical result for the law of the first intersecting line:
\begin{lemma}\label{loiprems}
Let $L$ denote the distance from $0$ to the first intersection on $\Delta_{0}$ of the line process,
and $\Theta$ the polar angle of this intersecting line, then the law of $(L,\Theta)$ is given by
$$\d(L,\Theta)(P)(\ell,\theta)={}{e^{-2\ell}}\,\cos\theta\,\un{\theta\in(-\pi/2,\pi/2)}
\,\d\ell\,\d\theta.$$
\end{lemma}
On the other hand, from formula \ref{formule-defaut-app} we get easily that
\begin{eqnarray*}
{d}_{\lambda}(0)&=& \overline{d}_{\lambda}(0)\\
&=&\lambda^{2}R(0,\omega)\cos\Theta(0,\omega)\sqrt{1+\frac{2L(0,\omega)\cos\Theta(0,\omega)}{\lambda^{2}R(0,\omega)}}\\
&&-\lambda^{2}R(0,\omega)\sqrt{1-\sin^2{\Theta(0,\omega)}
\left(1+\frac{2L(0,\omega)\cos\Theta(0,\omega)}{\lambda^{2}R(0,\omega)}\right)}\\
&&-L(0,\omega)
\end{eqnarray*}
if the inner term of the square root is non 
negative ({\it i.e.\/} $\lambda\geq\sqrt{({2L\sin^{2}\Theta})/({R\cos\Theta})}\,$), 
$+\infty$ otherwise. The asymptotic expansion
of those square roots gives easily
\begin{eqnarray*}
\overline{d}_{\lambda}(0)&=&-
\frac{L(0,\omega)^{2}}{2\lambda^2}\,\frac{\cos{2}\Theta(0,\omega)}{R(0,\omega)\cos\Theta(0,\omega)}+O(\lambda^{-4}).
\end{eqnarray*}
 
\subsubsection*{Step 2: Almost sure covering}
We conclude the proof of theorem \ref{une-direction} by stating the following lemma:
\begin{lemma}\label{pas-trop-peu}
As $\delta$, $r$, $\epsilon$ tend to zero and $s$ tends to $+\infty$, one has
$$\P(\Omega_{\delta,\epsilon,r,s})\to 1.$$
\end{lemma}
The proof of this lemma comes directly from the properties of the Poisson point process $\xxx$ and the asymptotic results on the
law of the inner and outer radii of the Crofton cell stated in section 2.

\cqfd\end{proof}
\begin{remark}
The almost sure convergence above will not be used for the convergence of the defect process, only the convergence in law of
the finite directional distributions is needed, however we shall state them almost surely.
\end{remark}

\subsection{Two and more directions}

For more directions we may state similar results, 
\begin{theorem}\label{cv-globale}
For all $0\leq t_{1}<\cdots<t_{n}<2\pi$, the finite dimensional random vector
$\lambda^2(d_{\lambda}(t_{1},\cdot),\ldots,d_{\lambda}(t_{n},\cdot))$ 
converges almost surely towards $(Z({t_{1},\cdot}),\ldots,Z({t_{n}},\cdot))$, where the law of this
random vector may be fully explicited.
\end{theorem}
The proof is essentially the same one as for one direction, only with more technical details.\\

\begin{figure}[h]
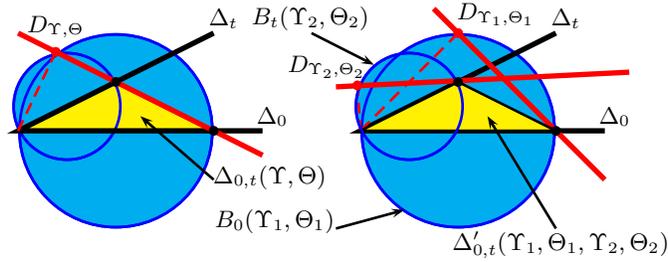

\caption{The sets $\Delta_{0,t}$, $\Delta'_{0,t}$, $B_{0}$ and $B_{t}$, and the corresponding
lines.}\label{biboule}
\end{figure}

This random vector depends only on the characteristics of the Crofton cell, let us for instance give the exact 
law of this vector for two directions (by invariance under rotations we choose directions $0$ and $t\in(0,2\pi)$):
$\lambda^2(d_\lambda(0,\cdot),d_\lambda(t,\cdot))$ converges in law
as $\lambda$ goes to infinity towards 
\begin{eqnarray*}
B\,\left(Z_{0}(\Upsilon,\Theta,R),Z_{t}(\Upsilon,\Theta,R)\right)
+(1-B)\left(W_{0}(\Upsilon_{1},\Theta_{1},R_{1}),
W_{t}(\Upsilon_{2},\Theta_{2},R_{2})\right),
\end{eqnarray*}
where 
\begin{itemize}
\item $B$ is Bernoulli random variable stating that the same line determines the intersections 
in directions
$0$ and $t$: this occurs with probability $p$,
\begin{eqnarray*}
p&=&\E\left[\sum_{(\rho,\alpha,R)\in\zzz}
\un{\text{The lines from }\zzz\setminus\{\rho,\alpha,R\}\text{ do not intersect $\Delta_{0}$ or 
$\Delta_{t}$ before }D_{(\rho,\alpha)}}\right],\\
&=&\int_{\RR_{+}\times[0,2\pi)\times[R_{\star},+\infty)}
\exp(-{\mathfrak p}(\Delta_{0,t}(\rho,\alpha)))\,\d\rho\,\d\alpha\,\d\mu(r),
\end{eqnarray*}
where $\Delta_{0,t}(\rho,\alpha)$ is the triangle described by figure \ref{biboule}, ${\mathfrak p}$ denotes the 
perimeter function. This Bernoulli random variable is independent from the following ones,
\item $(\Upsilon,\Theta,R)$ has the following distribution:
\begin{eqnarray*}
\d(\Upsilon,\Theta)(P)(\rho,\alpha,r)&=&p^{-1}\,\un{\alpha\in(-\pi/2,\pi/2),\alpha-t\in(-\pi/2,\pi/2)}\\
&&\exp(-{\mathfrak p}(\Delta_{0,t}(\rho,\alpha)))\,\d\rho\,\d\alpha\,\d\mu(r),
\end{eqnarray*}
\item $\displaystyle{Z_{0}(\Upsilon,\Theta,R)=-\frac{\Upsilon^{2}}{2R\cos^{2}\Theta}
\frac{\cos 2\Theta}{\cos\Theta}}$,
\item $\displaystyle{Z_{t}(\Upsilon,\Theta,R)=-\frac{\Upsilon^{2}}{2R\cos^{2}(\Theta-\theta)}
\frac{\cos2(\Theta-t)}{\cos(\Theta-t)}}$,
\item $(\Upsilon_{1},\Theta_{1},R_{1},\Upsilon_{2},\Theta_{2},R_{2})$ has the following distribution,
\begin{eqnarray*}
&&\d(\Upsilon_{1},\Theta_{1},R_{1},\Upsilon_{2},\Theta_{2},R_{2})(P)
(\rho_{1},\alpha_{1},r_{1},\rho_{2},\alpha_{2},r_{2})\\
&=&(1-p)^{-1}\,\un{\alpha_{1}\in(-\pi/2,\pi/2)}\un{\alpha_{2}-t\in(-\pi/2,\pi/2)}
\un{(\rho_{1},\alpha_{1})\notin B_{t}(\rho_{2},\alpha_{2})}\,
\un{(\rho_{2},\alpha_{2})\notin B_{0}(\rho_{1},\alpha_{1})}\\
&&\exp(-{\mathfrak p}(\Delta'_{0,t}(\rho_{1},\alpha_{1},\rho_{2},\alpha_{2})))\,
\d\rho_{1}\,\d\alpha_{1}\,\d\mu(r_{1})\,\d\rho_{2}\,\d\alpha_{2}\,\d\mu(r_{2}),
\end{eqnarray*}
where the sets $B_{0}$, $B_{t}$ and $\Delta'_{0,t}$ are described by figure \ref{biboule}.
\item $\displaystyle{W_{0}(\Upsilon_{1},\Theta_{1},R_{1})=-\frac{\Upsilon_{1}^{2}}{2R_{1}\cos^{2}\Theta_{1}}
\frac{\cos{2}\Theta_{1}}{\cos\Theta_{1}}}$,
\item $\displaystyle{W_{\theta}(\Upsilon_{2},\Theta_{2},R_{2})=-\frac{\Upsilon_{2}^{2}}{2R_{2}\cos^{2}
(\Theta_{2}-t)}
\frac{\cos{2}(\Theta_{2}-t)}{\cos(\Theta_{2}-t)}}$,
\end{itemize}

\vspace*{17mm}

\psset{unit=0.65cm}
\hspace*{2cm}\pscircle[linecolor=blue,fillcolor=cyan,fillstyle=solid,linewidth=1pt](2,0){2}
\pscircle[linecolor=blue,fillcolor=cyan,fillstyle=solid,linewidth=1pt](1,0.5){1.118}
\pscircle[linecolor=blue,fillcolor=cyan,fillstyle=none,linewidth=1pt](2,0){2}
\pscircle[linecolor=blue,fillcolor=cyan,fillstyle=solid,linewidth=1pt](9,0){2}
\pscircle[linecolor=blue,fillcolor=cyan,fillstyle=solid,linewidth=1pt](8,0.5){1.118}
\pscircle[linecolor=blue,fillcolor=cyan,fillstyle=none,linewidth=1pt](9,0){2}
\pspolygon[linecolor=black,linewidth=1pt,fillcolor=yellow,fillstyle=solid](2,1)(4,0)(0,0)
\pscircle[linecolor=blue,linewidth=1pt](1,0.5){1.118}
\psline[linecolor=black,linewidth=2pt](5,0)(0,0)(4,2)
\psline[linecolor=black,linewidth=2pt](12,0)(7,0)(11,2)
\psline[linecolor=red,linewidth=2pt](5,-0.5)(0,2)
\psline[linecolor=red,linewidth=1pt,linestyle=dashed](0,0)(0.77,1.615)
\pscircle[linecolor=red,linewidth=1pt,fillcolor=red,fillstyle=solid](0.77,1.615){0.1}
\psline[linecolor=black,linewidth=1pt]{->}(4,-1)(2.5,0.3)\rput[b](5.15,-1.2){\small $\Delta_{0,t}(\Upsilon,\Theta)$}
\rput[b](0.77,1.9){\small $D_{\Upsilon,\Theta}$}
\rput[b](5.2,0.1){\small $\Delta_0$}
\rput[b](4.2,2.1){\small $\Delta_{t}$}
\pscircle[linecolor=black,linewidth=1pt,fillcolor=black,fillstyle=solid](4,0){0.1}
\pscircle[linecolor=black,linewidth=1pt,fillcolor=black,fillstyle=solid](2,1){0.1}
\rput[b](12.2,0.1){\small $\Delta_0$}
\rput[b](11.2,2.1){\small $\Delta_{t}$}
\psline[linecolor=red,linewidth=2pt](6.5,0.9)(12.5,1.2)
\psline[linecolor=red,linewidth=2pt](12,-1)(8.5,2.5)
\psline[linecolor=red,linewidth=1pt,linestyle=dashed](7,0)(9,2)
\psline[linecolor=red,linewidth=1pt,linestyle=dashed](7,0)(6.94,0.93)
\pspolygon[linecolor=black,linewidth=1pt,fillcolor=yellow,fillstyle=solid](9,1)(11,0)(7,0)
\pscircle[linecolor=blue,linewidth=1pt](8,0.5){1.118}
\pscircle[linecolor=red,linewidth=1pt,fillcolor=red,fillstyle=solid](9,2){0.1}
\pscircle[linecolor=red,linewidth=1pt,fillcolor=red,fillstyle=solid](6.94,0.93){0.1}
\pscircle[linecolor=black,linewidth=1pt,fillcolor=black,fillstyle=solid](11,0){0.1}
\pscircle[linecolor=black,linewidth=1pt,fillcolor=black,fillstyle=solid](9,1){0.1}
\psline[linecolor=black,linewidth=1pt]{->}(11,-2)(9.6,0.3)\rput[b](11.1,-2.6){\small $\Delta'_{0,t}(\Upsilon_{1},\Theta_{1},\Upsilon_{2},\Theta_{2})$}
\rput[b](9.77,2.16){\small $D_{\Upsilon_{1},\Theta_{1}}$}
\rput[b](6.31,1.09){\small $D_{\Upsilon_{2},\Theta_{2}}$}
\psline[linecolor=black,linewidth=1pt]{->}(6.5,-2)(7.9,-1.62)\rput[b](5.25,-2.1){\small $B_{0}(\Upsilon_{1},\Theta_{1})$}
\psline[linecolor=black,linewidth=1pt]{->}(6,2)(7.28,1.35)\rput[b](5.98,2.05){\small $B_{t}(\Upsilon_{2},\Theta_{2})$}

\vspace*{1.1cm}

\section{Convergence of the stochastic process}

\begin{figure}[h]
\begin{center}
\begin{picture}(300.00,112.00)
\includegraphics{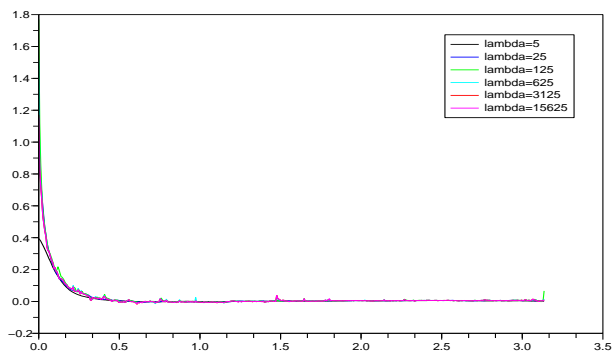}
\end{picture}
\end{center}
\vspace*{-4cm}
\caption{Covariograms, sample of size 250000}\label{covario}
\end{figure}

In this section we consider the processes 
$(\lambda^2d_{\lambda}(t))_{t\in[0,2\pi)}$. Let us first remark the following:
knowing the joint limit law of the couples $(\lambda^2d_{\lambda}(0),\lambda^2d_{\lambda}(t))$ gives some knowledge on this process,
for instance by simulation we can obtain the covariogram 
$t\mapsto \cov(\lambda^2d_{\lambda}(0),\lambda^2
d_{\lambda}(t))$, $t\in[0,\pi]$,
in figure \ref{covario}. One clearly observes the divergence as $\lambda$ tends to infinity
of the covariance for $t\to 0$, this is a consequence of the following elementary result coming from
the explicit law of the defect:
\begin{corollary}\label{integrabilite}
The limit expected defect is an integrable random variable with
$$\lim_{\lambda\to +\infty}
\E[\lambda^2 d_{\lambda}(0)]=0.$$
However, this limit expected defect is not square-integrable:
$$\lim_{\lambda\to+\infty}\E\left[\left(\lambda^2 d_{\lambda}(0)\right)^{2}\right]
=+\infty.$$
\end{corollary}

\vspace*{4.5cm}

The matter of convergence of the whole process
will be stated in the state $D([0,2\pi])$. As a matter of fact for each $\lambda$ the trajectory of the defect 
process is continous, however the limit process is not continuous: the choice of the space $D([0,2\pi])$,
even if there is no geometric justification in choosing right-continuity, seems to be quite natural.
Let us thus consider such processes $X_{\lambda}$ on $[0,2\pi]$, 
according to theorem 15.4 in \cite{bill} the conditions for the convergence of 
processes $(X_{\lambda})_{\lambda\geq 1}$ 
on $D([0,2\pi])$ are:
\begin{itemize}
\item convergence in law of the finite-dimensional distributions, this is true thanks to theorem \ref{cv-globale};
\item tightness criterion, for instance the following one: $\forall \eta,\epsilon>0$ there exists $\delta>0$ such that
$$\limsup_{\lambda\to+\infty} \P\left(
\sup_{t_{1}\leq t\leq t_{2},\ t_{2}-t_{1}<\delta}\min(
|X_{\lambda}(t)-X_{\lambda}(t_{1})|,|X_{\lambda}(t_{2})-X_{\lambda}(t)|)\geq \epsilon
\right)\leq \eta.$$
\end{itemize}

\begin{figure}[h]
\begin{tabular}{rl}
\begin{picture}(150.00,75.00)
\includegraphics{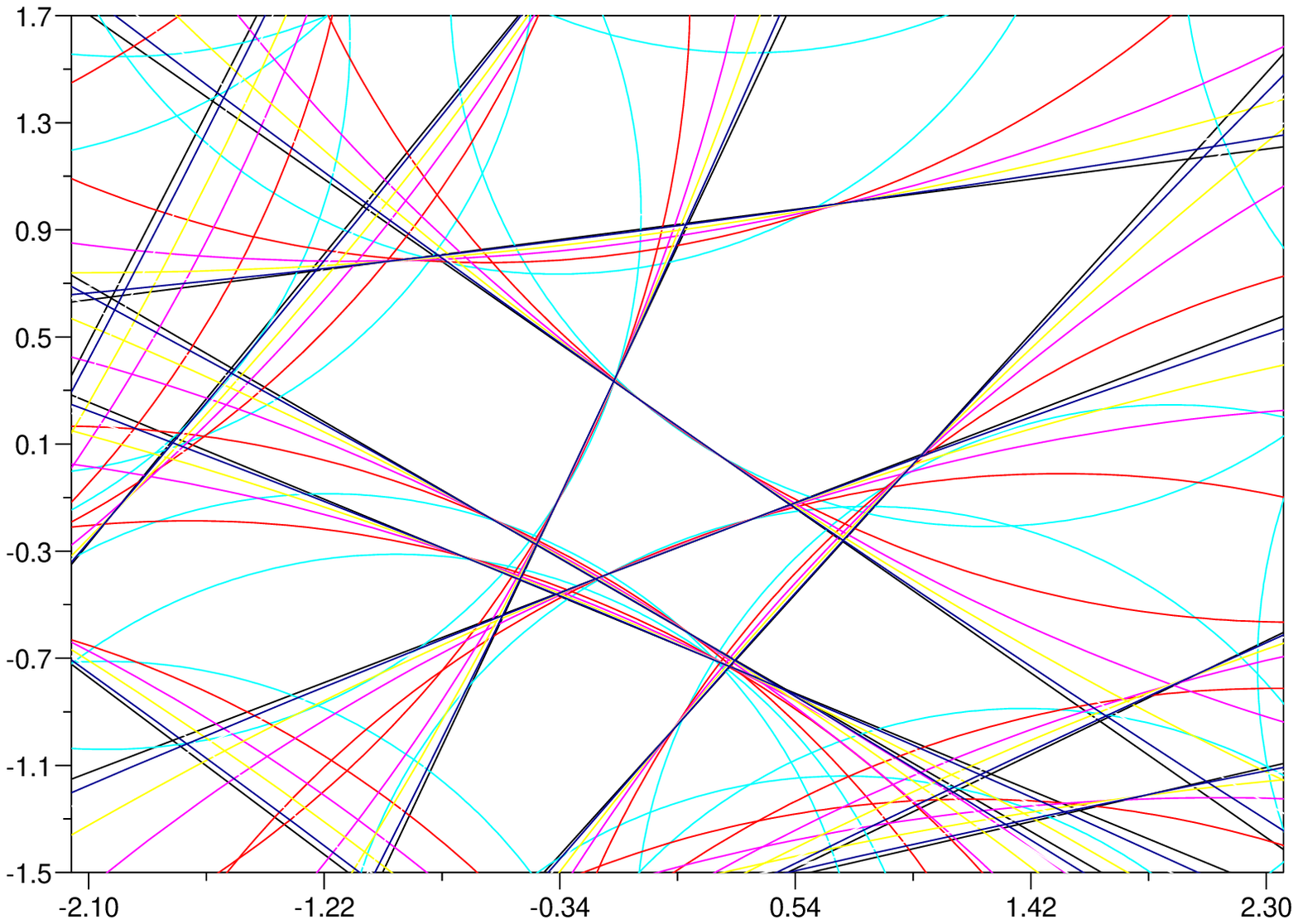}
\end{picture}
&
\begin{picture}(150,75)
\includegraphics{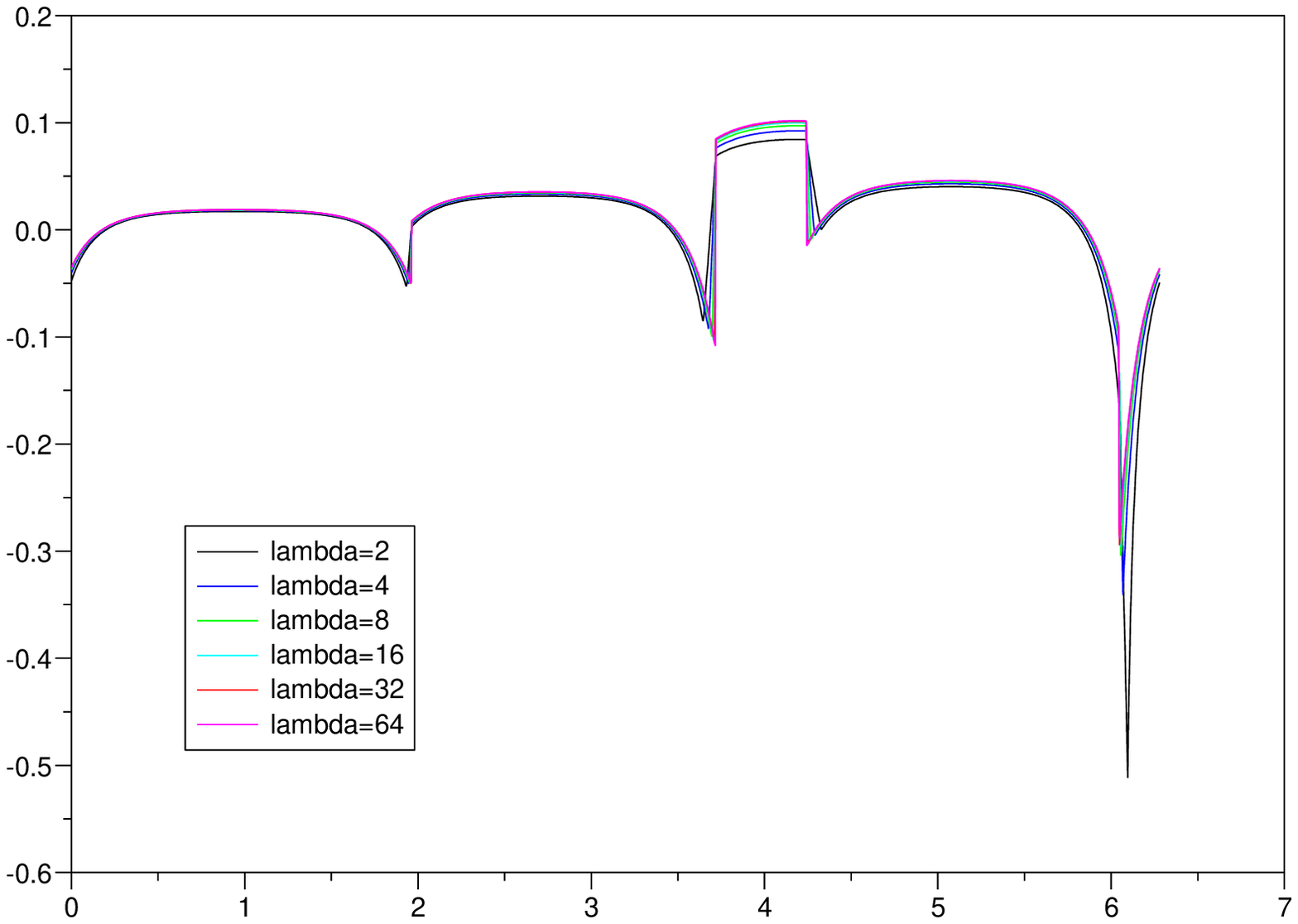}
\end{picture}
\end{tabular}
\vspace*{-3cm}
\caption{The Crofton cell and the related processes, one checks that there are four edges, and four
singularities}\label{melange}
\end{figure}

In our case, unfortunately one can not use directly such a tightness criterion: indeed if we take $X_{\lambda}=\lambda^{2}
d_{\lambda}$ we see (figure \ref{melange}) that the high slopes that appear near the angles corresponding to the 
vertices of the Crofton cell
forbid us to use this kind of citerion, as well as all other classical criteria. Hence we shall first
show the convergence of the approximate defect process $X_{\lambda}=\lambda^{2}\overline{d}_{\lambda}$, as this
process is the combination of a jump process and a smooth process, and then 
give an explicit estimate on the accuracy of this approximation in $L^1$ norm.

\vspace*{4cm}

\subsection{Convergence of the approximate defect process}

We prove the following theorem on the approximate defect:
\begin{theorem}\label{cv-def-app}
The approximate defect process $(\lambda^{2}\overline{d}_{\lambda}(t))_{t\in[0,2\pi]}$ 
converges in law in $D([0,2\pi])$ to the process
$(X_{t})_{t\in[0,2\pi]}$ defined for all $i\in\{1,\ldots,N\}$ and $t\in[\theta_{i},\theta_{i+1})$ by
$$\omega\mapsto X_{t}(\omega)=-\frac{\Upsilon_{i}(t,\omega)^{2}}{2R_{i}(t,\omega)}
\frac{\cos2(\Theta_{i}(t,\omega)-t)}{\cos^{3}(\Theta_{i}(t,\omega)-t)},$$
using the notations of definition \ref{def-def-app}.
\end{theorem}

\begin{proof}
Let us fix $\epsilon$ and $\eta$, both positive numbers, and define for $0<r<s$ and $\delta_{0}>0$ the
set $\Omega_{r,R,\delta_{0}}$ of those $\omega\in\Omega$ such that the Crofton
cell ${\mathscr C}(\omega)$ and $D^{\lambda}$ satisfy:
\begin{itemize}
\item $B(0,r)\subset {\mathscr C}(\omega)\subset B(0,s)$, {\it i.e.\/}
  $R_{m}\geq r$ and $R_{M}\leq s$;
\item $D^{\lambda}\subset B(0,s)$, {\it i.e.\/} $R_M(\lambda)\le s$;
\item the angular distance $\theta_{i+1}(\omega)-\theta_{i}(\omega)$ 
between any two consecutive vertices of ${\mathscr C}(\omega)$ is greater than $\delta_{0}$.
\end{itemize}
We give without proof the following lemma, similar to lemma \ref{pas-trop-peu}, 
stating that with high probability the Crofton cell is a `gentle' polygon:
\begin{lemma}\label{mauvais-petits}
As $r\to 0$, $s\to +\infty$ and $\delta_{0}\to 0$ one has $\P(\Omega_{r,s,\delta_{0}})\to 1$.
\end{lemma}

On this event $\Omega_{r,s,\delta_{0}}$ we check easily from definition \ref{def-def-app} that for $t\in[\theta_{i},\theta_{i+1})$
the approximate defect $\lambda^{2}\overline{d}_{\lambda}(t)$ is well defined
for $\lambda\geq \sqrt{2s^{3}/(r^{2}R_{\star})}$ and is Lipschitz-continuous for $\lambda\geq 2\sqrt{s^{3}/(r^{2}R_{\star})}$
with Lipschitz constant lesser than 
$13s^{6}/(r^{4}R_{\star})$
on the interval $[\theta_{i},\theta_{i+1})$.

It is then straightforward to check that if we choose $\delta<\delta_{0}$ and 
$\delta<C\epsilon r^{4}R_{\star}/s^{6}$, we have for $\lambda$ large enough (depending on $r$ and $s$)
\begin{eqnarray*}
&&\P\left(
\sup_{t_{1}\leq t\leq t_{2},\ t_{2}-t_{1}<\delta}\min(
|\lambda^{2}\overline{d}_{\lambda}(t)-\lambda^{2}\overline{d}_{\lambda}(t_{1})|,
|\lambda^{2}\overline{d}_{\lambda}(t_{2})-\lambda^{2}\overline{d}_{\lambda}(t)|)\geq \epsilon
\right)\\
&\leq &1-
\P(\Omega_{r,s,\delta_{0}}).
\end{eqnarray*}
We may then conclude using lemma \ref{mauvais-petits}.
\cqfd\end{proof}

The asymptotic expansion of the approximate defect gives also a convergence in the spaces $L^{p}(0,2\pi)$, $p\in [1,+\infty]$: indeed
on the event $\Omega_{r,s,\delta_{0}}$ one checks that  $(\lambda^{2}\overline{d}_{\lambda}(t)-X_{t})_{t\in(0,2\pi)}$
is bounded uniformely on $(0,2\pi)$ by $\lambda^{-2}s^{5}/(R_{\star}^{2}r^{2})$ times an explicit 
constant depending only on $r$, $s$ and $\delta_{0}$, hence
\begin{proposition}\label{cv-lp}
Almost surely one has for all $p\in[1,+\infty]$
$$(\lambda^{2}\overline{d}_{\lambda}(t))_{t\in(0,2\pi)}
\converge{$L^{p}$}{\lambda\to+\infty}(X_{t})_{t\in(0,2\pi)}.$$
\end{proposition}

\subsection{Estimate on the accuracy of the approximate defect process}\label{cv-lun-section}

Obviously the approximation of the rescaled defect process $(\lambda^{2}d_{\lambda}(t))_{t\in[0,2\pi]}$ 
by the process $(\lambda^{2}\overline{d}_{\lambda}(t))_{t\in[0,2\pi]}$ is not convergent to $0$ in the space $L^{\infty}(0,2\pi)$
because of the (common) jumps of both the rescaled and limit processes. We prove the following
result:
\begin{theorem}\label{cv-Lun} 
Almost surely one has the following convergence:
$$\lim_{\lambda\to+\infty}\lambda^{2}\int_{(0,2\pi)}|d_{\lambda}(t)-\overline{d}_{\lambda}(t)|\,dt=0.$$
\end{theorem}
\begin{proof}
The proof of this theorem will also be done in two steps:
\begin{itemize}
\item estimates on the widths of the `almost jumps' of the defect process on almost-full probability events;
\item $L^{\infty}$ estimates on the difference of the two processes on those same events.
\end{itemize}

\noindent{\em Step 1: widths of jumps}

The quantities $d_{\lambda}$ and $\overline{d}_{\lambda}$ are distinct only in the following case:
the first intersecting line in direction $t$ does not induce the first intersecting disc in direction $t$, this decomposes into three 
subcases:
\begin{itemize}
\item this first intersecting disc is associated to an other line of the line process that contains an edge of the 
Crofton cell, adjacent to the actual edge intersected by $\Delta_{t}$,
\item this disc is associated to a non-adjacent edge,
\item this disc is associated to a line that does not induce any edge of the Crofton cell.
\end{itemize} 
We shall show that the last two cases can be excluded on some event: let us first define for $\epsilon>0$ the 
thick Crofton cell as ${\mathscr C}_{\epsilon}={\mathscr C}\oplus B(0,\epsilon)$, we shall say that 
it is {\em equivalent} to ${\mathscr C}$ if it has the same edges and vertices as ${\mathscr C}$, more precisely if for each 
$i\in\{1,\ldots,N_{e}\}$ one has
$$\left(D_{i}\oplus B(0,\epsilon)\cap D_{i+1}\oplus B(0,\epsilon)\right)\cap
\left(D_{i+1}\oplus B(0,\epsilon)\cap D_{i+2}\oplus B(0,\epsilon)\right)=\emptyset,$$
where the $D_{i}$, $i\in\{1,\ldots,N_{e}\}$ are the lines supporting the edges of ${\mathscr C}$.
The intersections $D_{i}\oplus B(0,\epsilon)\cap D_{i+1}\oplus B(0,\epsilon)$ are lozenges, denoted by 
$C_{i,\epsilon}$ (crossings).
\begin{remark}
It is clear that by thickening the Poisson line process the thick Crofton cell is defined by at most the 
lines on the boundary of the 
Crofton cell. Our notion of equivalence is a little more demanding than just 
assuming that all those lines bound ${\mathscr C}_{\epsilon}$.
\end{remark}
\begin{lemma}\label{les-autres-mauvais}
Let $r<s$ and $\epsilon$ be positive numbers and define
$\overline{\Omega}_{r,s,\epsilon}$ the event such that for $\omega\in\overline{\Omega}_{r,s,\epsilon}$
\begin{itemize}
\item the Crofton cell ${\mathscr C}(\omega)$ is included in $B(0,s)$ and contains $B(0,r)$;
\item the thick crofton cell ${\mathscr C}_{\epsilon}(\omega)$ is equivalent to the Crofton cell ${\mathscr C}$, and
for each $(\rho,\theta,R)\in\zzz
\setminus\{(\Upsilon_{i},\Theta_{i},R_{i}),\ i\in\{1,\ldots,N_{e}\}\}$ one has
$${\mathscr C}_{\epsilon}(\omega)\cap \left(D_{\rho,\theta}\oplus B(0,\epsilon)\right)=\emptyset,$$
and for all $i\in\{1,\ldots,N_{e}\}$
$${C}_{i,\epsilon}\cap \left(D_{\rho,\theta}\oplus B(0,\epsilon)\right)=\emptyset,$$
\end{itemize}
Then one has
$$\lim_{r,\epsilon\to 0,s\to+\infty}\P(\overline{\Omega}_{r,s,\epsilon})=1.$$
Furthemore for $\lambda$ large enough (depending on $r$, $s$ and $\epsilon$ and $(\Upsilon_{i},\Theta_{i})$, $i\in
\{1,\ldots,N_{e}\}$) for all $t\in[0,2\pi)$
the first intersecting line and disc are associated either to the same point $(\Upsilon_{i},\Theta_{i},R_{i})\in \zzz$
or to the two points $(\Upsilon_{i},\Theta_{i},R_{i})$ and $(\Upsilon_{i\pm1},\Theta_{i\pm1},R_{i\pm1})$
where $i\in\{1,\ldots,N_{e}\}$ is the index of the corresponding edge of the Crofton cell.
\end{lemma}
The proof of this lemma follows classical lines, for instance for fixed $r$
and $s$ it is obvious that the conditional probability given that $R_{M}\leq s$ and $R_{m}\geq r$
that some line does violate the third or fourth hypothesis is
of order $\epsilon s$. The last point is clearly illustrated in
figure \ref{dessin-a-faire}.\\

If $\lambda$ is chosen large enough, depending only on $r$, $s$, $\epsilon$, and $(\Upsilon_{i},\Theta_{i})_{i\in\{1,\ldots,N_{e}\}}$, 
then it is clear thanks to proposition \ref{cv-hausdorff} that the angles $t$
for which the two edges are needed to determine the defect at $t$ are at most those corresponding to 
the disjoint sub-lozenges $C_{i,\eta/\lambda^{2}}$ formed by the intersections of two polar thick 
lines $D_{{i}}\oplus B(0,\eta/\lambda^{2})$
and $D_{i+1}\oplus B(0,\eta/\lambda^{2})$ corresponding to edges $i$ and $i+1$, where 
$\eta$ depends on $r$, $s$, $\epsilon$ and $(\Upsilon_{i},\Theta_{i})_{i\in\{1,\ldots,N_{e}\}}$, 
with $\eta/\lambda^{2}<\epsilon$.

\vspace*{4cm}

\psset{unit=1cm}
\newrgbcolor{grisclair}{0.8 0.8 0.8}
\newrgbcolor{grisfonce}{0.4 0.4 0.4 }
\psline[linewidth=4mm,linecolor=grisclair](2,2.7)(10,2.7)
\psline[linewidth=4mm,linecolor=grisclair](2,0)(6,4)
\psline[linewidth=4mm,linecolor=grisclair](10,0)(5,4)
\pspolygon[linewidth=1pt,linecolor=black,fillstyle=solid,fillcolor=grisfonce](4.2424,2.5)(4.7424,2.5)(5.1624,2.9)(4.6624,2.9)
\pspolygon[linewidth=1pt,linecolor=black,fillstyle=solid,fillcolor=grisfonce](6.54,2.5)(7.14,2.5)(6.66,2.9)(6.06,2.9)
\psline[linewidth=1pt,linecolor=black](10,0)(5,4)
\psline[linewidth=1pt,linecolor=black](2,0)(6,4)
\psline[linewidth=1pt,linecolor=black](2,2.7)(10,2.7)
\psline[linestyle=dashed,linewidth=1pt,linecolor=black](3.5,0)(6.06,2.9)
\psline[linestyle=dashed,linewidth=1pt,linecolor=black](3.5,0)(7.14,2.5)
\psline[linestyle=dashed,linewidth=1pt,linecolor=black](3.5,0)(4.2424,2.5)
\psline[linestyle=dashed,linewidth=1pt,linecolor=black](3.5,0)(5.1624,2.9)
\pscircle[linecolor=black,fillcolor=black,fillstyle=solid](3.5,0){0.05}
\rput[b](3.1,0){$0$}
\begin{figure}[h]
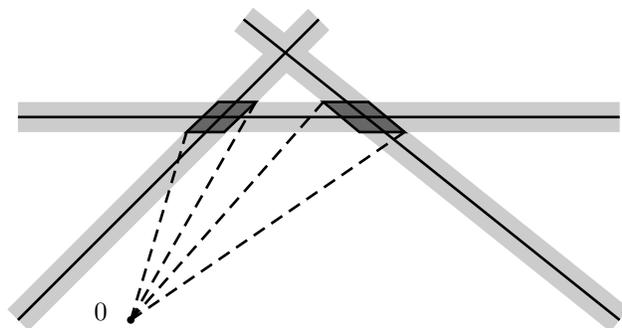

\caption{Thick Crofton cell and lozenges.}\label{dessin-a-faire}
\end{figure}

This implies that the total length of problematic angles $t$ is lesser than $N\times \zeta/\lambda^{2}$, where $\zeta$ depends 
in quite a technical way on $r$, $s$, $\epsilon$ and ${\mathscr C}$ through the minimum of the 
differences $|\Theta_{i+1}-\Theta_{i}|$. 

\noindent{\em Step 2: $L^{\infty}$ bounds}
From the properties of the Crofton cell on the set $\overline{\Omega}_{r,s,\epsilon}$, we may easily evaluate exactly the
difference between both processes: indeed when this difference is non zero it is lesser than
\begin{eqnarray*}
|d_{\lambda}(t)-\overline{d}_{\lambda}(t)|&\leq&
d_{1}(t)+d_{2}(t)+d_{3}(t),
\end{eqnarray*}
where $d_{1}(t)$ (resp. $d_{2}(t)$) is the modulus of the distance between the intersections on $\Delta_{t}$ of
$D_{i}$ and the associated disc $B_{i}$ (resp. $D_{i\pm1}$ and $B_{i\pm1}$), and $d_{3}(t)$ is the modulus of the
distance between $D_{i}\cap\Delta_{t}$ and $D_{i\pm1}\cap\Delta_{t}$ (see figure \ref{dessin-a-faire} above).\\

Let us remark that one of $d_{1}(t)$ and $d_{2}(t)$ is exactly
$\overline{d}_{\lambda}(t)$: thanks to proposition \ref{cv-hausdorff}, those two terms are bounded from above by a constant for $\lambda$ large enough:
$$d_{j}(t)\leq
\frac{1}{\lambda^{2}R_*}\frac{s}{r}\left(s+\frac{s^2}{\lambda^2R_*}\right)^2, \quad j=1,2.$$
The third term is also bounded from above by a constant $M(r,s,\epsilon,{\mathscr C})\times \eta/\lambda^{2}$.
Hence we obtain
$$\lambda^{2}|d_{\lambda}(t)-\overline{d}_{\lambda}(t)|\leq M_{1}(r,s,\epsilon,{\mathscr C}).$$

Hence we have the following estimate on the $L^{1}$ norm of the difference:
\begin{eqnarray*}
\lambda^{2}\int_{(0,2\pi)}|d_{\lambda}(t)-\overline{d}_{\lambda}(t)|\,dt &\leq&\frac{1}{\lambda^{2}}\,
M_{2}(r,s,\epsilon,{\mathscr C})\times{N_{e}},
\end{eqnarray*}
this upper bound converges towards zero for each $\omega$ in the set 
$\overline{\Omega}_{r,s,\epsilon}$: this concludes the proof of theorem \ref{cv-Lun}.
\cqfd\end{proof}

As a consequence of this result and of proposition \ref{cv-lp} we obtain eventually the following convergence
\begin{theorem}
Almost surely one has
$$(\lambda^{2}d_{\lambda}(t))_{t\in(0,2\pi)}\converge{$L^{1}(0,2\pi)$}{\lambda\to+\infty}(X_{t})_{t\in(0,2\pi)},$$
where the process $(X_{t})_{t\in(0,2\pi)}$ is defined in theorem \ref{cv-def-app}.
\end{theorem}

\subsection{Tail probability for the supremum of the defect process}

This short section is devoted to a uniform bound on the tail probabilities for
the defect processes for large $\lambda$'s:
\begin{proposition}\label{prop-sup}
One has the following estimate: 
for all $\beta<1$ 
$$\lim_{u\to+\infty} \left\{s^{\beta}\limsup_{\lambda\to+\infty}
\P\left(\sup_{\theta\in[0,2\pi)}\lambda^2|d_{\lambda}(\theta)|\geq u\right)\right\}=0.$$
\end{proposition}
\begin{proof}[Sketch]
The proof uses the same tools as before, the estimates on the growth of both
${\mathscr C}$ and the empty component of the Boolean model around the origin. 
Indeed we may give an explicit upper bound of the defect
on the set $\overline{\Omega}_{r,s,\epsilon}$ using the computations of section \ref{cv-lun-section}. 
This bound is roughly of order $s^{5}/r^{2}$. 
By using an adequate choice of $r\to 0$ and $s\to+\infty$
in terms of powers of $u$ we may obtain the result. This tedious proof is left to the reader.
\cqfd
\end{proof}

\begin{remark}[Splitting the defect process]
The limit defect process may be decomposed in a continuous part and a pure jump part, such a splitting can be
done for the defect process at fixed $\lambda$: let us indeed write
$$\lambda^{2}d_{\lambda}(t)=\lambda^{2}\overline{\overline{d}}_{\lambda}(t)+
\left(\lambda^{2}d_{\lambda}(t)-\lambda^{2}\overline{\overline{d}}_{\lambda}(t)\right),$$
where $\overline{\overline{d}}_{\lambda}$ is the continuous part of the process $\overline{d}_{\lambda}$
(this is almost surely defined as being equal to $\overline{d}_{\lambda}$ at the angle $t=0$, and the jumps are {\em deleted}).
Then it can be shown that both terms above converge in law, the first one in the space $C(0,2\pi)$, and the second
one in a weak sense.
\end{remark}

\begin{remark}[Directions for the general case]
For more general shapes the coupling is more tricky to obtain, we may proceed in the following way:
\begin{itemize}
\item consider smooth shapes with no flat portion on the boundary: 
$\tau(F)$, where $\tau$ is a uniform rotation, and $F$ a smooth random closed set;
\item given $\tau$ and $F$, assign to a Poisson line the centre of the rotated rescaled shape tangent at the line, at the same distance from
the origin than the line;
\item compute the intensity of the point process of centres of shapes, and modify this intensity so that it becomes 
the Lebesgue measure multiplied by the parameter $\lambda^{2}$, this shall be done by a function
$$(\rho,\theta,\tau,F)\mapsto \Psi_{\lambda}(\rho,\theta,\tau,F)\in{\mathbb R}_{+}\times[0,2\pi);$$
\end{itemize}
Using this procedure, the computations might be done involving more technical details, the
limiting process might be expressed with the curvature of $F$.
\end{remark}

\bibliographystyle{apt}

\bibliography{nouvelle-corrigee-270409}

\end{document}